\documentclass[12pt]{article}
\usepackage{amsmath}
\usepackage{amssymb}
\usepackage{tabularx}
\usepackage{enumerate}
\usepackage{graphicx}

\newtheorem{thm}{Theorem}[section]
\newtheorem{lem}[thm]{Lemma}
\newtheorem{cor}[thm]{Corollary}
\newtheorem{prop}[thm]{Proposition}
\newtheorem{rmk}{Remark}[section]
\newtheorem{defi}{Definition}[section]
\newtheorem{exm}{Example}[section]
\newtheorem{pppp}{Proof}

\newcommand{\qed}{\hspace{1em}\mbox{\raisebox{0.65ex}{\fbox{}}}}

\numberwithin{equation}{section}

\newcommand{\be}{\begin{equation}}
\newcommand{\ee}{\end{equation}}
\newcommand\bes{\begin{eqnarray}} \newcommand\ees{\end{eqnarray}}
\newcommand{\bess}{\begin{eqnarray*}}
\newcommand{\eess}{\end{eqnarray*}}

\newcommand{\R}{\mathbb{R}}
\newcommand{\bpf}{{\bf Proof:\ \ }}
\newcommand{\epf}{\mbox{}\hfill $\Box$}

\begin{document}

\thispagestyle{empty}

\title{A SIS reaction-diffusion-advection model in a low-risk and high-risk domain
\thanks{The work is partially supported by the
NRF of Korea (Grant No. 2010-0025700), the NSFC of China (Grant No. 11371311), the High-End Talent Plan of Yangzhou University, China and NSERC of Canada.}}
\date{\empty}

\author{Jing Ge$^{a}$, Kwang Ik Kim$^{b}$, Zhigui Lin$^a$ and Huaiping Zhu$^{c}$\\
{\small $^a$School of Mathematical Science, Yangzhou University, Yangzhou 225002, China}\\
{\small $^b$Department of Mathematics, Pohang University of Science and Technology, } \\
{\small Pohang, 790-784, Republic of Korea}\\
{\small $^c$Laboratory of Mathematical Parallel Systems (LAMPS),}\\ 
 {\small Department of Mathematics and Statistics,}\\
 {\small York University, Toronto, ON, M3J 1P3, Canada}
}

 \maketitle

\begin{quote}
\noindent
{\bf Abstract.} { 
A simplified SIS model is proposed and investigated to understand the impact of spatial heterogeneity of environment and advection on the persistence and eradication of an infectious disease. The free boundary is introduced to model the spreading front of the disease. The basic reproduction number associated with the diseases in the spatial setting is introduced. Sufficient conditions for the disease to be eradicated or to spread are given. Our result shows that if the spreading domain is high-risk at some time, the disease will continue to spread till the whole area is infected; while if the spreading domain is low-risk, the disease may be vanishing or keep spreading depends on the expanding capability and the initial number of the infective individuals. The spreading speeds are also given when spreading happens, numerical simulations are presented to illustrate the impacts of the advection and the expanding capability on the spreading fronts.
 }


\end{quote}

\section{Introduction}

Mathematical models have been made to investigate the
transmission of infectious diseases and the asymptotic profiles of the steady states of the diseases
(see \cite{AM1, AM3, ZH}).  For classical
compartmental epidemic models for infectious diseases described by ordinary differential systems,
it is well-known 
that the so called basic reproduction number determines 
whether the disease will be endemic \cite{BC, vanden}.
It is also common that for vector-borne diseases,
backward bifurcation may occur in the compartmental models which reveals that besides the basic reproduction number, the endemic also
depends on the initial sizes of the involving individuals ( see \cite{Wenzhang98, WZhu} and references therein).
 In recent years, spatial diffusion and environmental heterogeneity have been recognized as important factors to affect 
 the persistence and eradication 
 of infectious diseases such as measles, tuberculosis and flu,  etc.,
 especially for vector-borne diseases, such as malaria, dengue fever, West Nile virus etc.
More importantly, it is the spatial transmission and environmental heterogeneity that decide the speed and pattern of the spatial spread of infectious diseases.
In this case, the usual basic reproduction number will not be enough to describe the disease transmission dynamics,  especially to reflect the spatial features of the spread in the region considered.  Therefore, it is essential to investigate the role of diffusion on the transmission
and control of diseases in a heterogeneous environment.

To understand the dynamics of disease transmission in
a spatially heterogeneous environment, an SIS epidemic reaction-diffusion model has
been proposed by Allen, Bolker, Lou and Nevai in \cite{AL}, and the model is described by the following coupled parabolic system:
\begin{eqnarray}
\left\{
\begin{array}{lll}
S_{t}-d_S\Delta S=-\frac{\beta (x) SI}{S+I}+\gamma(x)I,\; &  x\in\Omega, \ t>0, \\
I_{t}-d_I\Delta I=\frac{\beta (x) SI}{S+I}-\gamma(x)I,\; &   x\in \Omega,\ t>0
\end{array} \right.
\label{Aa1}
\end{eqnarray}
with homogeneous Neumann boundary condition
$$\frac {\partial S}{\partial \eta}=\frac {\partial I}{\partial \eta}=0,\  x\in \partial \Omega, \ t>0,$$
where $S(x, t)$ and $I(x, t)$  represent the density of susceptible and infected individuals at
location $x$ and time $t$, respectively, the positive constants $d_S$ and $d_I$ denote the corresponding diffusion rates for
the susceptible and infected populations, $\beta(x)$ and $\gamma(x)$ are positive H$\ddot{o}$lder continuous functions,
 which account for spatial dependent rates of disease contact transmission and disease recovery at $x$, respectively.
 The term $\frac{\beta (x) SI}{S+I}$ is the standard incidence of disease, since the term $\frac{S I}{S + I}$ is a
 Lipschitz continuous function of $S$ and $I$ in
the open first quadrant, it can be extended to define to the entire first quadrant by defining it to be
zero when either $S = 0$ or $I = 0$.

Letting $N=S+I$, adding two equations in (\ref{Aa1}) and then integrating over $\Omega$ yields
that $\frac{\partial}{\partial t}\int_\Omega N(x,t)dx=0$ for $t>0$, this means that the total population size remains a constant,
 and recovered individuals 
 become susceptible after survival from the infectious of the disease.

 As in \cite{AL},
we say that $x$ is a {\bf high-risk site} if the local disease transmission rate $\beta(x)$ is greater
than the local disease recovery rate $\gamma(x)$. An {\bf low-risk site} is defined in a similar manner.
The habitat $\Omega$ is characterized as {\bf high-risk}
( or {\bf low-risk} ) if the spatial average $(\frac{1}{|\Omega|}\int_{\Omega}\beta(x)dx)$
of the transmission rate is greater than ( or less than ) the spatial
average $(\frac{1}{|\Omega|}\int_{\Omega}\gamma(x)dx)$ of the recovery rate, respectively.

To characterize the dynamics of the transmission of the disease, the authors in \cite{AL} introduced the basic reproduction number
$R^N_0$ (we use $R^N_0$ for Neumann boundary condition to compare with $R^D_0$ for Dirichlet boundary condition defined later) by
\begin{equation}\label{AAq}R^N_0=R^N_0(\Omega)=\ \sup_{\phi\in H^1(\Omega),\, \phi\neq 0}\{\frac{\int_\Omega \beta \phi^2 dx}{\int_\Omega (d_I|\triangledown \phi|^2+\gamma \phi^2)dx}\}.\end{equation}
They showed that if $R_0^N < 1$, the population density
$(S(t,\cdot), I(t,\cdot))$ converges to a unique disease free equilibrium $(S_0, 0)$, while
if $R_0^N > 1$, there exists a unique positive endemic equilibrium $(S^*, I^*)$. In addition,
 global stability of endemic steady state for some particular cases
and particularly the asymptotical profiles of the endemic steady states as the diffusion coefficient for susceptible
individuals is sufficiently small are given.
In some recent work \cite{Pe1, Pe3, PZ}, Peng et al. further investigated the asymptotic
behavior and global stability of the endemic equilibrium for system (\ref{Aa1}) subject to the Neumann boundary conditions,
and provided much understanding of the impacts of large and small diffusion rates of the susceptible
and infected population on the persistence and extinction of the
disease.

For the SIS reaction-diffusion model  (\ref{Aa1}) with Dirichlet boundary conditions
$$S=I=0,\quad  x\in \partial \Omega, \ t>0,$$
adding the two equations in (\ref{Aa1}), and integrating over $\Omega$ one yields
that $\frac{\partial}{\partial t}\int_\Omega N(x,t)dx \leq 0$ for $t>0$, it follows that the total population is decreasing and
$\int_\Omega N(x,t)dx\to 0$ as $t\to \infty$. To avoid the loss of population in the boundary and the diffusion process,
 Huang et al. \cite{HH} added an additional growth term in the first equation of (\ref{Aa1}) and studied the global
dynamics of the corresponding problem.

To focus on the new phenomena induced by spatial heterogeneity of environment, we assume that the population $N(x,t)$ is constant in space for all time, that is, $N(x,t)\equiv N^*$ for $x\in \Omega$ and $t\geq 0$, then the system (\ref{Aa1}) becomes the simplified SIS diffusive equation
\begin{equation}
I_{t}-d_I \Delta I=(\beta (x) -\gamma(x))I-\frac {\beta (x)}{N^*} I^2,\;  x\in \Omega,\ t>0.
\label{Aa2}
\end{equation}
As we know the solutions of equation (\ref{Aa2}) subject to the Neumann or Dirichlet boundary conditions
are always positive for any $t>0$  no matter what the nonnegative nontrivial initial data is given.
Epidemiologically, it means that the disease spreads and becomes endemic to the whole area immediately even the infection is limited in a small region at the beginning. It does not reflect the reality that disease always spreads gradually from an endemic region to spread further to an larger area in terms of spatial spread.

To describe the changing process of the domain, the free boundary problems
have been discussed in many areas (\cite{CS}), especially
 the well-known Stefan condition has been used to describe the interaction and spreading
 process at the boundary. For example, it was used
to describe the melting of ice in contact with water \cite{R}, the oxygen in muscles in \cite{C}, the wound healing in \cite{CF}, the spreading of
 the invasion species in \cite{DG, DGP, DL, GW, KY, LIN, MYY, Wa, WZ, ZX}.
 Recently, it was used to study an SIR epidemic model in a homogeneous environment in \cite{KLZ}.

For emerging and re-emerging infectious diseases, usually the spread of the disease starts at a source location and spread over areas where
contact transmission occurs.  For example, West Nile virus, a kind of mosquito-borne virus arrived and caused an encephalitis outbreak in New York city in 1999 \cite{CDCwnv99}. With mosquito as vector and bird as amplification host,  West Nile virus has kept spreading and become established in the North America continent \cite{CDCwnv}.  For infectious diseases like West Nile virus, it is essential to understand how the disease is spreading spatially
over further to larger area to cause endemic, to determine the condition for the virus to spread spatially, to predict the spatial spread for the purpose of prevention and control.

Normally, diffusion of particles in physics is random and obeys Fick's law.
However, species in population dynamics or diseases in epidemiology diffuse differently owing to their initiative behaviors and activities.
Some species or diseases prefer to move towards one direction because of appropriate climate, wind direction, etc.
For example, in studying the propagation of West
Nile Virus in North America, it was observed in \cite{Maiyang} that West Nile Virus appeared for the first time in
New York City in the summer of 1999. In the second year the wave front travels 187km to the north and 1100km to the south, till 2002,
it has been spread across almost the whole America continent.
Therefore, the propagation of WNv from New York City to California state is a consequence of the diffusion and advection movements of birds.
Especially, bird advection becomes an important factor for lower mosquito biting rates.

As one preliminary study, we will focus on the changing of the infected domain and the advection movement of the disease, and consider an SIS epidemic model
with the free boundary to describe the spreading frontier of the disease. Spatial advective diffusion and environmental heterogeneity are two very complex aspect of the spread of the infectious diseases.  For simplicity, we assume the region or environment is one dimensional. We will investigate the behavior of the positive solution
$(I(x, t); g(t), h(t))$ to the following problem
\begin{eqnarray}
\left\{
\begin{array}{lll}
I_{t}-d_I I_{xx}+\alpha I_{x}=(\beta (x) -\gamma(x))I-\frac {\beta (x)}{N^*} I^2,\; &g(t)<x<h(t),\; t>0,   \\
I(g(t),t)=0,\, g'(t)=-\mu I_{x}(g(t),t),&t>0,\\
I(h(t),t)=0,\, h'(t)=-\mu I_{x}(h(t),t),  & t>0, \\
g(0)=-h_0,\, h(0)=h_0,\, I(x,0)=I_{0}(x), & -h_0\leq x\leq h_0,
\end{array} \right.
\label{a3}
\end{eqnarray}
where $x=g(t)$ and $x=h(t)$ are the moving left and right
boundaries to be defined,  $h_0,\, d,\, \alpha$ and $\mu $ are positive constants.
$\alpha$ and $\mu$ are referred as the advection rate and the expanding capability, respectively.
As above, $\beta(r)$ and $\gamma(r)$
 are positive H$\ddot{o}$lder continuous functions, which account for spatial dependent rates of disease contact transmission and disease recovery at $x$, respectively.
 Further, in the paper we assume
$$(H)\;\;\;\;\lim_{x\rightarrow\pm \infty}\beta(x)=\beta_{\infty},\;\;\;
\lim_{x\rightarrow \pm \infty}\gamma(x)=\gamma_{\infty}\ \textrm{and}\ \beta_{\infty}-\gamma_{\infty}>0,$$
which means that far sites of the habitat are similar and high-risk.

 In this paper, we only consider the small advection and assume that $\alpha<2\sqrt{(\beta_\infty-\gamma_\infty)d_I}$, it is well known that $2\sqrt{(\beta_\infty-\gamma_\infty)d_I}$ is the
minimal speed of the traveling waves to the cauchy problem
\begin{equation}
\label{tw11}
 u_t-d_I u_{xx}=u(\beta_\infty-\gamma_\infty-bu),\; t>0,\; x\in \R.
\end{equation}

The initial distribution of the infected class $I_0(x)$ is nonnegative and satisfies
\begin{equation}
I_0\in C^2[-h_0, h_0],\, I_0(-h_0)=I_0(h_0)=0\, \textrm{and} \ 0<I_0(x)\leq N^*,\, x\in (-h_0, h_0),
\label{Ae2}
\end{equation}
where the condition (\ref{Ae2}) indicates that at the beginning, the infected exists in the area with $x\in (-h_0, h_0)$,
but for the area $|x|\geq h_0$, no infected happens yet. Therefore,
the model means that beyond the left boundary $x=g(t)$ and the right boundary $x=h(t)$, there is only susceptible, no infectious individuals.

The equation governing the free boundary, the spread front,
$h'(t)=-\mu u_{r}(h(t),t)$, is a special case of the well-known
Stefan condition, which has been established in \cite{LIN} for the diffusive populations.
The positive constant $\mu$ measures the ability of the infected transmit and diffuse towards the new area.

Different from the usual compartmental models and reaction-diffusion models with Dirichlet boundary conditions,
it is natural that the basic reproduction number for the disease transmission modeled by the free boundary conditions will be time dependent.
For the reaction-diffusion models with free boundary conditions, we will define the basic reproduction number based on the definition
for Dirichlet boundary conditions, and use the basic reproduction number to characterize the dynamics of the temporal and spatial spread of the disease.
As a preliminary study, we will consider the case when the domain is one-dimensional and heterogenous, and focus to describe when
the diseases can be vanishing (eradicated) or can spread to become endemic further over the domain.

This paper is organized as follows. In the next section, the global
existence and uniqueness of the solution to
(\ref{a3}) are presented by using a contraction mapping theorem,
comparison principle is also employed. Section 3 is devoted to developing the basic reproduction
numbers and their properties. Sufficient conditions for the disease to vanish is given in section 4. Section 5 deals with
the case and conditions for the disease to spread and become endemic. Section 6 is devoted to the asymptotic spreading speed.
Numerical simulations are also given in section 7
to illustrate the impacts of the advection and the expanding capability on the free boundary, and a brief discussion is also presented.

\section{Existence and uniqueness}

In this section, we first prove the following local existence and
uniqueness result by the contraction mapping theorem. We then use
suitable estimates to show that the solution is well defined for all
$t>0$.

\begin{thm} For any given $I_0$ satisfying \eqref{Ae2},  and any $\nu \in (0, 1)$, there is a $T>0$ such that  problem \eqref{a3}
admits a unique solution
$$(I; g, h)\in C^{ 1+\nu,(1+\nu)/2}([g(t), h(t)]\times [0,T])\times C^{1+\nu/2}([0,T])\times C^{1+\nu/2}([0,T]);$$
moreover,
\begin{eqnarray} \|I\|_{C^{1+\nu,
(1+\nu)/2}([g(t), h(t)]\times [0,T])}+\|g\|_{C^{1+\nu/2}([0,T])}+\|h\|_{C^{1+\nu/2}([0,T])}\leq
C,\label{b12}
\end{eqnarray}
where $C$ and $T$ only depend on $h_0, \nu$ and $\|I_0\|_{C^{2}([-h_0, h_0])}$.
\end{thm}
\bpf  As in \cite{CF}, we first straighten the free boundary. Let
$\xi (y)$ be a function in $C^3[0, \infty)$ satisfying
$$\xi (y)=1\ \  \textrm{if} \,\, |y-h_0|<\frac{h_0}8,\; \xi (y)=0\ \
 \textrm{if} \,\, |y-h_0|>\frac {h_0}2,\quad |\xi
'(y)|<\frac 5{h_0} \mbox{ for all } y.$$
Consider the transformation $(y, t)\rightarrow (x, t)$, where
$$x=y+\xi (y)(h(t)-h_0)+\xi (-y)(g(t)+h_0),
\quad -\infty<y<+\infty.$$
As long as $|h(t)-h_0|\leq \frac {h_0}8$ and $|g(t)+h_0|\leq \frac {h_0}8$,
the above transformation $x \to y$ is a diffeomorphism from $(-\infty, +\infty)$
onto $(-\infty, +\infty)$. Moreover, it changes the left free boundary $x=g(t)$
to the line $x=-h_0$ and the right free boundary $x=h(t)$
to the line $x=h_0$. It follows from direct calculations that
\begin{eqnarray*}
\displaystyle\frac {\partial y}{\partial x}=\frac
1{1+\xi'(y)(h(t)-h_0)-\xi'(-y)(g(t)+h_0)}&\equiv &
A(g(t),h(t), y),\\
\displaystyle \frac {\partial^2 y}{\partial x^2}=-\frac
{\xi''(y)(h(t)-h_0)+\xi'(-y)(g(t)+h_0)}{[1+\xi'(y)(h(t)-h_0)-\xi'(-y)(g(t)+h_0)]^3}&\equiv &
B(g(t), h(t), y),\\
\displaystyle \frac {\partial y}{\partial t}=\frac
{-\xi(y)h'(t)-\xi(-y)g'(t)}{1+\xi'(y)(h(t)-h_0)-\xi'(-y)(g(t)+h_0)}&\equiv &C(g(t), h(t), y).
\end{eqnarray*}

If we set
$$I(x, t)=y+\xi (y)(h(t)-h_0)+\xi (-y)(g(t)+h_0)=v(y, t),$$
then the free boundary problem \eqref{a3} becomes
\begin{eqnarray}
\left\{
\begin{array}{lll}
v_{t}-A^2d_I v_{yy}-(Bd_I-C-\alpha A)v_y=(\beta -\gamma)v-\frac {\beta }{N^*} v^2,\; &
-h_0<y<h_0,\; t>0, \\
v(g(t),t)=0,\, g'(t)=-\mu v_y(-h_0, t),&t>0, \\
v(h(t),t)=0,\, h'(t)=-\mu v_y(h_0, t),&t>0, \\
g(0)=-h_0,\, h(0)=h_0,\, v(y, 0)=v_0(y), &-h_0\leq y\leq h_0,
\end{array} \right.
\label{Bb}
\end{eqnarray}
where $A=A(g(t), h(t), y)$, $B=B(g(t), h(t), y)$, $C=C(g(t), h(t), y)$ and $v_0=I_0$.

The rest of the proof follows from
the contraction mapping theorem together with standard $L^p$ theory and the Sobolev imbedding theorem
\cite{LSU}, we then omit it here, see Theorem 2.1 in \cite{DL} for details. \epf

To show that the local solution obtained in Theorem 2.1 can be
extended to all $t>0$, we need the following estimate.

\begin{lem} Let $(I; g, h)$ be a solution to problem \eqref{a3} defined for $t\in (0,T_0]$ for some $T_0>0$.
Then we have
\[
0<I(x, t)\leq N^*\; \mbox{ for } g(t)<x<
h(t),\; t\in (0,T_0]. \]
\end{lem}
\bpf
It is easy to see that $I\geq 0$ in $[g(t), h(t)]\times
[0, T_0]$ as long as the solution exists. Using the strong maximum principle to the equations in $[g(t), h(t)]\times
[0, T_0]$ yields
 $$I(x, t)>0\;\; \textrm{for} \ g(t)<x<h(t),\, 0< t\leq T_0.$$

Since $I(x,t)$ satisfies
\begin{eqnarray*}
\left\{
\begin{array}{lll}
I_{t}-d_I I_{xx}+\alpha I_x\leq \beta (r) I(1-\frac {I}{N^*}), &
g(t)<x<h(t),\, t>0,  \\
I(g(t), t)=I(h(t),t)=0,\quad &t>0\\
I(r,0)=I_0(r)\leq N^*,\; &-h_0\leq r\leq h_0,
\end{array} \right.
\end{eqnarray*}
direct application of the maximum principle gives that $I\leq N^*$ in $[g(t), h(t)]\times [0, T_0]$.
\epf

The next lemma shows that the left free boundary for problem (\ref{a3}) is strictly monotone decreasing and the right free boundary is increasing.

\begin{lem} \label{mono}Let $(I; g, h)$ be a solution to problem \eqref{a3} defined for $t\in (0,T_0]$ for some $T_0>0$.
Then there exists a constant $C_1$  independent of $T_0$ such
that
\[
-C_1\leq g'(t)<0\; \mbox{and}\; 0<h'(t)\leq C_1 \; \mbox{ for } \; t\in (0,T_0]. \]
\end{lem}
\bpf
Using the strong maximum principle to the equation of $I$ gives that
 $$I_x(g(t),t)>0 \, \mbox{and}\, I_x(h(t), t)<0 \ \;\; \textrm{for} \ 0<t\leq T_0.$$
 Hence $g'(t)<0$ and $h'(t)>0$ for $t\in (0, T_0]$ by using the free boundary condition in (\ref{a3}).

It remains to show that $g'(t)\geq C_1$ and $h'(t)\leq C_1$ for $t\in (0,T_0]$ and
some $C_1$. The proof is similar as that of Lemma 2.2 in \cite{DL} with $C_1=2MN^*\mu$ and
$$M=\max\left\{\frac{\alpha}{d_I}+\sqrt{\frac{\overline \beta}{2d_I}},\
\frac{4\|I_0\|_{C^1([-h_0,h_0])}}{3N^*}\right\},\quad \overline \beta =\max_{[-h_0, h_0]} \beta(r),$$ we omit it here.
 \epf

\begin{thm} The solution of problem \eqref{a3} exists and is
unique for all $t\in (0,\infty)$.
\end{thm}
\bpf
It follows from the uniqueness and Zorn's lemma that there is a number $T_{\max}$
such that $[0, T_{\max})$ is the maximal time interval in which the
solution exists. Now we prove that $T_{\max}=\infty$ by the contradiction argument. Assume that
$T_{\max}<\infty$. By Lemmas 2.2 and 2.3, there exist $C_1$ independent of $T_{\max}$ such that for
$x\in [g(t), h(t)]$ and $t\in [0, T_{\max})$ ,
\[ 0\leq I(x,t)\leq N^*, \;(x,t)\in [g(t), h(t)]\times [0, T_{\max}),\]
\[ -h_0-C_1t\leq g(t)\leq -h_0, \ -C_1\leq g'(t)\leq 0,\; t\in [0, T_{\max}),\]
\[ h_0\leq h(t)\leq h_0+C_1 t, \ 0\leq h'(t)\leq C_1,\; t\in [0, T_{\max}).\]
We now fix $\delta_0\in (0, T_{\max})$ and $M>T_{\max}$. By standard
parabolic regularity, we can find $C_2>0$ depending only on
$\delta_0$, $M$ and $C_1$ such that
 $$||I(\cdot, t)||_{C^2[g(t), h(t)]}\leq
C_2$$
 for $t\in [\delta_0, T_{\max})$.  It then follows from the
proof of Theorem 2.1 that there exists a $\tau>0$ depending only on
$C_1$ and $C_2$ such that the solution of  problem \eqref{a3}
with initial time $T_{\max}-\tau/2$ can be extended uniquely to the
time $T_{\max}-\tau /2+\tau$. But this contradicts the assumption.
The proof is complete. \epf

In what follows, we shall exhibit the comparison principle, which is similar to Lemma 3.5 in \cite{DL}.
\begin{lem} (The Comparison Principle)
  Assume that $T\in (0,\infty)$, $\overline g, \overline
h, \underline g, \underline
h$ $\in C^1([0,T])$, $\overline I(x,t)\in C([\overline g(t), \overline h(t)]\times [0, T])\cap
C^{2,1}((\overline g(t), \overline h(t))\times (0, T])$, $\underline I(x,t)\in C([\underline g(t), \underline h(t)]\times [0, T])\cap
C^{2,1}((\underline g(t), \underline h(t))\times (0, T])$, and
\begin{eqnarray*}
\left\{
\begin{array}{lll}
\overline I_{t}-d_I \overline I_{xx}+\alpha \overline I_x\geq (\beta (x) -\gamma(x))\overline I-\frac {\beta (x)}{N^*} \overline I^2,
&\overline g(t)<x<\overline h(t), \ 0<t\leq T,\\
\underline I_{t}-d_I \underline I_{xx}+\alpha \underline I_x\leq (\beta (x) -\gamma(x))\underline I-\frac {\beta (x)}{N^*} \underline I^2,
&\underline g(t)<x<\underline h(t), \ 0<t\leq T,\\
\overline I(\overline g(t), t)=0,\; \overline g'(t)\leq -\mu \overline I_x(\overline g(t), t),\quad & 0<t\leq T,\\
\underline I(\underline g(t), t)=0,\; \underline g'(t)\geq -\mu \underline I_x(\underline g(t), t),\quad & 0<t\leq T,\\
\overline I(\overline h(t), t)=0,\; \overline h'(t)\geq -\mu \overline I_x(\overline h(t), t),\quad & 0<t\leq T,\\
\underline I(\underline h(t), t)=0,\; \underline h'(t)\leq -\mu \underline I_x(\underline h(t), t),\quad & 0<t\leq T,\\
\overline g(0)\leq -h_0<h_0\leq \overline h(0), \, I_{0}(x)\leq \overline I(x, 0),\;  &-h_0\leq x\leq h_0,\\
-h_0\leq \underline g(0)\leq \underline h(0)\leq h_0, \, \underline I(x,0)\leq I_{0}(x),\;  &\underline g(0)\leq x\leq \overline h(0),
\end{array} \right.
\end{eqnarray*}
then the solution $(I(x,t); g(t), h(t))$ to the free boundary problem $(\ref{a3})$ satisfies
$$\overline g(t)\leq g(t)\leq\underline g(t),\ \underline h(t)\leq h(t)\leq\overline h(t)\, \mbox{in}\, (0, T],$$
$$\underline I(x, t) \leq I(x, t)\, \mbox{for}\,  (x, t)\in [\underline g(t), \underline h(t)]\times (0, T],$$
$$I(x, t) \leq \overline I(x, t)\, \mbox{for}\,  (x, t)\in [g(t), h(t)]\times (0, T].$$
\end{lem}

\bigskip
The pair $(\overline u; \overline g, \overline h)$ in Lemma 2.5 is usually called an upper solution
of the problem \eqref{a3} and $(\underline u; \underline g, \underline h)$ is then called a lower solution.
To examine the dependence of the solution on the expanding capability $\mu$,
 we write the solution as $(I^{\mu}; g^{\mu}, h^{\mu})$. As a corollary of Lemma 2.5, we have the following monotonicity:

\begin{cor} For fixed $I_0, \alpha, h_0, \beta (x)$ and $\gamma (x)$.
If $\mu_1\leq \mu_2$. Then $I^{\mu_1}(x, t)\leq I^{\mu_2}(x, t)$ in $[g^{\mu_1}(t), h^{\mu_1}(t)]\times (0, \infty)$
 and $g^{\mu_2}(t)\leq g^{\mu_1}(t)$, $h^{\mu_1}(t)\leq h^{\mu_2}(t)$ in $(0, \infty)$.
\end{cor}

\section{Basic reproduction numbers}

In this section, we first present the basic reproduction number and its properties and implications
for the reaction-diffusion system \eqref{Aa2} with Dirichlet boundary condition, and then discuss the basic reproduction number
for the free boundary problem \eqref{a3}.

Let us introduced the basic reproduction number $R_0^D$ by
$$R_0^D=R_0^D(\Omega, d_I)=\ \sup_{\phi\in H^1_0(\Omega),\phi\neq 0} \{\frac{\int_{\Omega} \beta
\phi^2dx}{\int_{\Omega} (d_I|\triangledown \phi|^2+\gamma \phi^2)dx}\},$$
 the following result was given in \cite{HH} (Lemma 2.3):

\begin{lem} $1-R_0^D$ has the same sign as $\lambda_0$, where $\lambda_0$ is the principle eigenvalue of the reaction-diffusion problem
\begin{eqnarray}
\left\{
\begin{array}{lll}
-d_I\Delta \psi=\beta(x) \psi-\gamma(x)\psi+\lambda \psi,\; &
x\in \Omega,  \\
\psi(x)=0, &x\in \partial \Omega.
\end{array} \right.
\label{B1f}
\end{eqnarray}
\end{lem}

With the above defined reproduction number, we have

\begin{thm} The following assertions hold.

$(a)$ $R_0^D$ is a positive and monotone decreasing function of $d_I$;

$(b)$ $R_0^D\to \max_{x\in \Omega}\frac {\beta (x)}{\gamma (x)}$ as $d_I\to 0$;

$(c)$ $R_0^D\to 0$ as $d_I\to \infty$;

$(d)$ There exists a threshold value $d^*_I\in [0, \infty)$ such that $R_0^D>1$ for $d_I<d_I^*$ and $R_0<1$ for $d_I>d_I^*$.
 If all sites in the domain are lower-risk $(\beta(x)\leq \gamma(x)$ for $x\in \Omega)$, we have $R_0^D<1$ for all $d_I>0$;

$(e)$ Let $B_h$ be a ball in $R^n$ with the radius $h$. Then $R_0^D(B_{h})$ is strictly monotone increasing function of $h$, that is if $h_1<h_2$, then $R_0(B_{h_1})<R_0(B_{h_2})$.
Moreover, $\lim_{h\to \infty} R_0^D(B_{h})\geq \frac {\beta_{\infty}}{\gamma_{\infty}}$ provided that $(H)$ holds;

$(f)$ If $\Omega =(-h_0,h_0)$, $\beta (x)\equiv \beta^*$ and $\gamma(x)\equiv \gamma^*$, then
$$R_0^D=\frac { \beta^*}{d_I(\frac {\pi}{2h_0})^2+\gamma^* }.$$
\label{basic1}
\end{thm}
\bpf
The proof of part (a), (b) and (d) are similar to that of Theorem 2 in \cite{AL}. The threshold value in part (d) can be described
in the following manner:
$$d^*_I=\ \sup \{\frac{\int_{\Omega} (\beta-\gamma)
\phi^2dx}{\int_{\Omega} |\triangledown \phi|^2dx}:\ \phi\in W^{1,2}_0(\Omega)\ \textrm{and}\, \int_{\Omega} (\beta-\gamma)
\phi^2dx>0\}.$$
It is easy to see that if $\beta(x)\leq \gamma(x)$ for $x\in \Omega$, then $d_I^*=0$.

Next, let's first established part (f).
It is well-known fact that the principle eigenvalue of the problem
\begin{eqnarray*}
\left\{
\begin{array}{lll}
-d_I\Delta \psi=\lambda \psi,\; &
x\in (-h_0, h_0),  \\
\psi(x)=0, &x=\pm h_0
\end{array} \right.
\end{eqnarray*}
is $d_I(\frac {\pi}{2h_0})^2 $, the desired result follows if $\beta $ and $\gamma$ are constants.

The proof of the monotonicity in (e) is similar to that of Corollary 2.3 in \cite{CC}.
For the limit in (e), it follows the assumption $(H)$ that for any $\varepsilon>0$,
there exists a positive constant $r_0$ such that for $|x|\geq r_0$,
$$|\beta(x)-\beta_\infty|<\varepsilon,\ |\gamma(x)-\gamma_\infty|<\varepsilon.$$
Let $\phi_h(r)$ be the function in $C^2[0,h]$ satisfying
$$\phi_h(r)=1\, \textrm{if}\, |r|\leq h-\frac 34,\ \phi_h(r)=0\, \textrm{if}\, |r|\geq h-\frac 14,\ |\phi_h'(r)\leq 4\, \textrm{for}\, r\in [h-1,h].$$
By the definition of $R^D_0$,  we have
\begin{eqnarray*}
R_0^D(B_h)&\geq&  \frac{\int_{B_h} \beta(x)
\phi^2_h(|x|)dx}{\int_{B_h} (d_I|\triangledown \phi_h|^2+\gamma(x)\phi^2_h(|x|))dx}\\
&=& \frac{\int_{B_h/B_{h-1}} \beta(x)
\phi^2_h(|x|)dx+(\int_{B_{h-1}/B_{r_0}}+\int_{B_{r_0}}) \beta(x)dx}{\int_{B_h/B_{h-1}} d_I|\triangledown \phi_h|^2 dx
+(\int_{B_{h}/B_{r_0}}+\int_{B_{r_0}}) \gamma(x)dx}\\
&\geq &\frac{(\beta_\infty-\varepsilon)|B_{h-1}/B_{r_0}|}{4d_I|B_h/B_{h-1}|
+(\gamma_\infty+\varepsilon)|B_h/B_{r_0}|+\max_{x\in B_{r_0}}|B_{r_0}|}
\end{eqnarray*}
for $h>r_0$, therefore
\begin{eqnarray*} \liminf_{h\to +\infty} \ R_0^D(B_h)&\geq& \liminf_{h\to +\infty} \ \frac{(\beta_\infty-\varepsilon)|B_{h-1}/B_{r_0}|}{4d_I|B_h/B_{h-1}|
+(\gamma_\infty+\varepsilon)|B_h/B_{r_0}|+\max_{x\in B_{r_0}}|B_{r_0}|}\\
&=&\frac{\beta_\infty-\varepsilon}{\gamma_\infty+\varepsilon},
\end{eqnarray*}
which together with the monotonicity of $R_0^D(B_h)$ and the arbitrariness of small $\varepsilon$ gives  $\lim_{h\to \infty}\, R_0^D(B_{h})\geq \frac {\beta_{\infty}}{\gamma_{\infty}}$.

It remains to established part (c).  Now we show that $R_0^D\to 0$ as $d_I\to \infty$. In fact, if it is not true,
there exists a positive $a>0$ such that  $R_0^D\geq a$ for any $d_I>0$ since $R_0^D$ is monotone decreasing function of $d$.
It is a well-known fact that there exists a positive function $\phi(x)\in C^2(\Omega)$ such that $||\phi||_{L^\infty}=1$ and
\begin{eqnarray*}
\left\{
\begin{array}{lll}
-d_I\Delta \phi+\gamma \phi=\frac {\beta}{R_0^D} \phi,\; & x\in \Omega,  \\
\phi(x)=0, &x\in \partial \Omega.
\end{array} \right.
\end{eqnarray*}
Dividing both sides of the above equation by $d_I$ yields
$$ -\Delta \phi+\frac{\gamma}{d_I} \phi=\frac {\beta}{R_0^Dd_I} \phi.$$
Since $\frac{\gamma}{d_I}\to 0$ and $\frac {\beta}{R_0^Dd_I}\to 0$ as $d_I\to \infty$,
it follows from elliptic regularity that $\phi\to \overline \phi$ in $C(\Omega)$ as $d\to \infty$ for some positive function
$\overline \phi$ satisfying
$$ -\Delta \overline \phi=0 \ \textrm{in}\ \Omega, \quad \overline \phi=0 \ \textrm{on} \ \partial \Omega,$$
we then have $\overline \phi \equiv 0$ in $\Omega$, which leads to a contradiction.
\epf

\bigskip
For the following reaction-diffusion-advection problem
\begin{eqnarray}
\left\{
\begin{array}{lll}
I_t-d_I I_{xx}+\alpha I_x=(\beta(x) -\gamma(x))I-\frac{\beta(x)}{N^*}I^2, &
x\in (g_0, h_0),\, t>0,  \\
I(x,t)=0, &x=g_0\, \textrm{or}\, x=h_0, t>0,
\end{array} \right.
\label{Aa2ad}
\end{eqnarray}
where $g_0<h_0$, let us introduced the basic reproduction number $R_0^{DA}$ by
$$R_0^{DA}=R_0^{DA}((g_0, h_0), d_I, \alpha)=\ \sup_{\phi\in H^1_0(g_0, h_0),\phi\neq 0} \{\frac{\int_{g_0}^{h_0} \beta
e^{\alpha x/d_I}\phi^2dx}{\int_{g_0}^{h_0} (d_Ie^{\alpha x/d_I}\phi_x^2+\gamma e^{\alpha x/d_I}\phi^2)dx}\}.$$
If $\phi\in H^1_0(g_0, h_0)$, then $\psi=e^{\alpha x/(2d_I)}\phi\in H^1_0(g_0, h_0)$ also, and the mapping $\phi \mapsto e^{\alpha x/(2d_I)}\phi$ is bijective,
therefore taking $\phi=e^{-\alpha x/(2d_I)}\psi$ gives that
$$R_0^{DA}=R_0^{DA}((g_0, h_0), d_I, \alpha)=\ \sup_{\psi\in H^1_0(g_0, h_0),\psi\neq 0} \{\frac{\int_{g_0}^{h_0} \beta
\psi^2dx}{\int_{g_0}^{h_0} (d_I\psi_x^2+\frac{\alpha^2}{4d_I}\psi^2+\gamma \psi^2)dx}\}.$$

The following result is from variational methods, see for example, Chapter 2 in \cite{CC}.
\begin{lem} $1-R_0^{DA}$ has the same sign as $\lambda_0$, where $\lambda_0$ is the principle eigenvalue of the reaction-diffusion-advection problem
\begin{eqnarray}
\left\{
\begin{array}{lll}
-d_I \psi_{xx}+\alpha \psi_x=\beta(x) \psi-\gamma(x)\psi+\lambda \psi,\; &
x\in (g_0, h_0),  \\
\psi(x)=0, &x=g_0\, \textrm{or}\, x=h_0.
\end{array} \right.
\label{B11f}
\end{eqnarray}
\end{lem}

Combining Theorem \ref{basic1} with the above defined reproduction number yields
\begin{thm} The following assertions hold.

$(a)$ $R_0^{DA}$ is a positive and monotone decreasing function of $\alpha$;

$(b)$ If $\alpha\neq 0$, $R_0^{DA}\to 0$ as $d_I\to 0$ or as $d_I\to \infty$;

$(c)$ If $\Omega_1\subseteqq \Omega_2 \subseteqq R^1$, then $R_0^{DA}(\Omega_1)\leq R_0^{DA}(\Omega_2)$, with strict inequality if $\Omega_2 \setminus \Omega_1$ is an open set. Moreover, $\lim_{(h_0-g_0)\to \infty}\, R_0^{DA}((g_0, h_0), d_I,\alpha)\geq \frac {\beta_{\infty}}{\frac{\alpha^2}{4d_I^2}+\gamma_{\infty}}$ provided that $(H)$ holds;

$(d)$ If $\beta (x)\equiv \beta_{\infty}$ and $\gamma(x)\equiv \gamma_{\infty}$, then
$$R_0^{DA}=\frac { \beta_{\infty}}{d_I(\frac {\pi}{h_0-g_0})^2+\frac{\alpha^2}{4d_I}+\gamma_{\infty} }.$$
\label{basic2}
\end{thm}

Noticing that the domain for the free boundary problem \eqref{a3} is changing with $t$, so the basic reproduction number is not a constant and
should be changing.
Now we introduced the basic reproduction number $R_0^F(t)$ for the free boundary problem \eqref{a3} by
$$R_0^F(t):=R_0^{DA}((g(t),h(t)),d_I,\alpha)=\ \sup_{\psi\in H^1_0(g(t),h(t)),\psi\neq 0} \{\frac{\int_{g(t)}^{h(t)} \beta
\psi^2dx}{\int_{g(t)}^{h(t)} (d_I\psi_x^2+\frac{\alpha^2}{4d_I}\psi^2+\gamma \psi^2)dx}\},$$
it follows from Lemma \ref{mono} and Theorem \ref{basic2} that
\begin{thm} $R_0^F(t)$ is strictly monotone increasing function of $t$, that is if $t_1<t_2$, then $R_0^F(t_1)<R_0^F(t_2)$.
Moreover, if $(H)$ holds and $h(t)-g(t)\to \infty$ as $t\to \infty$, then $\lim_{t\to \infty}\, R_0^F(t)\geq \frac {\beta _{\infty}}{\frac{\alpha^2}{4d_I^2}+\gamma _{\infty}}$.
\end{thm}

\begin{rmk} In this paper, we have assumed that $(H)$ holds and $\alpha<2\sqrt{(\beta_\infty-\gamma_\infty)d_I}$. By Theorem 3.5, we have
 $R_0^F(t_0)>1$ for some $t_0>0$ provided that $h(t)-g(t)\to \infty$ as $t\to \infty$.\label{rem1}
\end{rmk}

\section{Disease vanishing}
It follows from Lemma \ref{mono} that $x=g(t)$ is monotonic decreasing and $x=h(t)$ is monotonic increasing, so
there exist $g_\infty\in [-\infty, -h_0)$ and $h_\infty\in (h_0, +\infty]$ such that $\lim_{t\to +\infty} \ g(t)=g_\infty$
and $\lim_{t\to +\infty} \ h(t)=h_\infty$. The next lemma shows that both $h_\infty$ and
$g_\infty$ are finite or infinite simultaneously.
\begin{lem}
If $h_\infty<\infty$ or $g_\infty>-\infty$, then both $h_\infty$ and
$g_\infty$ are finite and
$$
R_0^{DA}((g_\infty, h_\infty),d_I, \alpha)\leq 1\ \textrm{and}\
\lim_{t\to\infty} \|I(\cdot, t)\|_{C([g(t),\,h(t)])}=0.$$
\label{spr}
\end{lem}
\bpf Without loss of generality, we assume that $h_\infty<\infty$, and prove that $R_0^{DA}\leq 1$, which implies that $g_\infty>-\infty$ by Remark \ref{rem1}.

Assume that $R_0^{DA}((g_\infty, h_\infty), d_I, \alpha)>1$ by contradiction. To see the dependence of $R^{DA}_0$ on the recovery rate $\gamma(x)$, 
we write $R_0^{DA}$ as $R_0^{DA}((g_\infty, h_\infty), \gamma(x))$. It follows from the continuity that there exists $T^*\gg 1$ such that
$R_0^{DA}((g(T^*), h(T^*)), \gamma(x))> 1$. Furthermore, for small $\varepsilon$,  $R_0^{DA}((g(T^*), h(T^*)), \gamma(x)+\varepsilon)> 1$, where 
$\varepsilon$ depends on $T^*$. 
Let $w(x,t)$ be the solution of
\begin{eqnarray}
\left\{
\begin{array}{lll}
w_{t}-d_I w_{xx}+\alpha w_x=w(\beta (x)-\gamma(x)-\varepsilon-\frac {\beta(x)}{N^*}w), &
g(T*)<x<h(T^*),\, t>T^*,  \\
w(g(T^*), t)=w(h(T^*),t)=0,\quad &t>T^*\\
w(x,T^*)=I(x, T^*),\; &g(T^*)\leq x\leq g(T^*),
\end{array} \right.\label{dfg}
\end{eqnarray}
direct application of the maximum principle gives that $I(x,t)\geq e^ {\varepsilon (t-T^*)}w(x,t)$ in $[g(T^*), h(T^*)]\times [T^*, \infty)$.

On the other hand, since that $R_0^{DA}((g(T^*), h(T^*)), \gamma(x)+\varepsilon)> 1$, by the method of upper and lower solutions and
its associated monotone iterations \cite{Pao}, we have $\lim_{t\to \infty} w(x,t)\to w_s(x)$ uniformly on $[g(T^*), h(T^*)$, where $w_s$ is the unique positive steady-state solution of problem ({\ref{dfg}) and satisfies
\begin{eqnarray*}
\left\{
\begin{array}{lll}
-d_I w''_{s}+\alpha w'_s=w_s(\beta (x)-\gamma(x)-\varepsilon-\frac {\beta(x)}{N^*}w_s), &
g(T*)<x<h(T^*), \\
w_s(g(T^*))=I(h(T^*))=0,\quad &
\end{array} \right.
\end{eqnarray*}
Therefore $\lim_{t\to \infty} w(0,t)=w_s(0)>0$, which together with $I(0,t)\geq e^ {\varepsilon (t-T^*)}w(0,t)$ gives that $\lim_{t\to \infty} I(0,t)=\infty$. This contradicts with the fact that $I\leq N^*$.

{\bf Step 2.} $\lim_{t\to
+\infty} \ ||I(\cdot,t)||_{C([g(t),\, h(t)])}=0$.

Let $\overline I(x,t)$ denote the unique solution of the problem
\begin{eqnarray}
\left\{
\begin{array}{lll}
\overline I_t -d_I \overline I_{xx}+\alpha \overline I_x=\overline I(\beta(x)-\gamma(x))-\frac {\beta(x)}{N^*}\overline I^2,\;
&g_\infty<x<h_\infty,\, t>0, \\
\overline I(g_\infty,0)=0, \quad \overline I(h_\infty,0)=0,  & t>0,\\
\overline I(x,0)=\tilde I_0(x), &g_\infty<x<h_\infty,
\end{array} \right.
\label{m23}
\end{eqnarray}
with
\begin{eqnarray*}
\tilde I_0(x)= \left\{
\begin{array}{lll}
I_0(x)&g_0\leq x\leq h_0, \\
0, & \mbox{ otherwise}.
\end{array} \right.
\end{eqnarray*}
The comparison principle gives $0\leq I(t,x)\leq \overline I(t, x)$
for $t>0$ and $x\in [g(t), h(t)]$.

Using the fact $R_0^{DA}((g_\infty, h_\infty),d_I, \alpha)\leq 1$, we find that
$0$ is the unique nonnegative steady-state solution of the problem (\ref{m23}). Choosing the lower solution as $0$ and upper solution as $\max\{ ||\tilde I_0(x)||_{L^\infty[g_\infty, h_\infty]},\, N^*\}$, it is shown, by the method of upper and lower solutions and
its associated monotone iterations, that the time-dependent solution converges to the unique nonnegative steady-state solution.
Therefore, $\overline u(x,t)\to 0$ uniformly
for $x\in [g_\infty, h_\infty]$ as $t\to \infty$ and then $\lim_{t\to
+\infty} \ ||I(\cdot,t)||_{C([g(t),\, h(t)])}=0$.
 \epf

Therefore the spatial transmission of a disease depends on whether  $h_\infty-g_\infty=\infty$ and $\lim_{t\to
+\infty} \ ||I(\cdot, t)||_{C([g(t), h(t)])}=0$. We then have the following definitions:

\begin{defi}
The disease is {\bf vanishing} if $h_\infty -g_\infty<\infty$ and
 $\lim_{t\to +\infty} \ ||I(\cdot, t)||_{C([g(t),h(t)])}=0$, while the disease is {\bf spreading} if $h_\infty -g_\infty=\infty$ and
 $\limsup_{t\to +\infty} \ ||I(\cdot, t)||_{C([g(t),h(t)])}>0$.
\end{defi}

 The next result shows that if $h_\infty-g_\infty<\infty$, then vanishing happens.
\begin{lem}  If $h_\infty-g_\infty<\infty$, then $\lim_{t\to
+\infty} \ ||I(\cdot, t)||_{C([0, h(t)])}=0$.
\end{lem}
\bpf
 This result can be proved by the same argument as Lemma 4.1 in \cite{LLZ} with obvious modification,
we omit it here for brevity.
\epf

Now we give sufficient conditions so that the disease is vanishing.
\begin{thm} Suppose $R_0^F(0)(:=R_0^{DA}((-h_0, h_0),d_I, \alpha))<1$.  Then $h_\infty-g_\infty<\infty$ and
$$\lim_{t\to +\infty} \ ||I(\cdot, t)||_{C([g(t), h(t)])}=0$$
 if $||I_0(x)||_{C([-h_0, h_0])}$ is sufficiently small.
\end{thm}
\bpf The proof is constructing a suitable upper solution for $I$, which is similar to that of Lemma 5.3 in \cite{LLZ}.
We give the sketch here for completeness.

Since that $R_0^{DA}((-h_0, h_0), d_I, \alpha)<1$, it follows from the continuity that there exists $\delta_0>0$ such that $R_0^{DA}((-h_0, h_0),d_I, \rho))\leq [R_0^{DA}((-h_0, h_0), d_I, \alpha)+1]/2<1$ for $|\rho-\alpha|\leq \delta_0$, therefore using Lemma 3.1 gives that there is a $\lambda_0>0$ and $\psi(x)>0$ in $(-h_0, h_0)$ such that
\begin{eqnarray}
\left\{
\begin{array}{lll}
-d_I\Delta \psi+\rho \psi=\beta(r) \psi-\gamma(r)\psi+\lambda_0 \psi,\; &
-h_0<x<h_0,  \\
\psi(x)=0, &x=\pm h_0.
\end{array} \right.
\label{B1f1}
\end{eqnarray}
Therefore, there exists a small $\delta >0$ such that
$$\delta (1+\delta)^2+[(1+\delta)^2-1]\overline \beta\leq \lambda_0,$$
where $\overline \beta=||\beta(r)||_{L^\infty[0,\infty)}$.

Similarly as in Lemma 3.8 in \cite{DL}, we set
$$\sigma (t)=h_0(1+\delta-\frac \delta 2 e^{-\delta t}), \  t\geq 0,$$
and $$w(t, x)=\varepsilon e^{-\delta t}\psi(rh_0/\sigma (t)), \ 0\leq
r\leq \sigma(t),\ t\geq 0.$$
We can choose $\varepsilon$ sufficient small such that if $||I_0||_{L^\infty}\leq \varepsilon \psi(\frac {h_0}{1+\delta/2})$, then
$(w(x,t), -\sigma(t), \sigma(t))$ be an upper solution of problem (\ref{a3}).
Applying Lemma 2.5 gives that $g(t)\geq-\sigma(t)$, $h(t)\leq\sigma(t)$
and $I(x, t)\leq w(x, t)$ for $g(t)\leq x\leq h(t)$ and $t>0$. It
follows that $h_\infty\leq \lim_{t\to\infty}
\sigma(t)=h_0(1+\delta)<\infty$, $g_\infty\geq -\sigma(t)>-\infty$ and then $\lim_{t\to
+\infty} \ ||I(\cdot, t)||_{C([g(t), h(t)])}=0$.
 \epf

From the above proof, we have the following result, see also Lemma 3.8 in \cite{DL} or Lemma 2.9 in \cite{DG}.
\begin{thm} Suppose $R_0^F(0)(:=R_0^{DA}((-h_0, h_0),d_I, \alpha))<1$.  Then $h_\infty-g_\infty<\infty$ and
$$\lim_{t\to +\infty} \ ||I(\cdot, t)||_{C([g(t), h(t)])}=0$$
 if $\mu$ is sufficiently small.
\end{thm}

\section{Disease spreading}

In this section, we are going to give the sufficient conditions so that the disease is spreading. We first
prove that if $R_0^F(0)(:=R_0^{DA}((-h_0, h_0),d_I, \alpha))\geq 1$, the disease is spreading.
\begin{thm} If $R_0^F(0)\geq 1$, then $h_\infty=-g_\infty=\infty$ and $\liminf_{t\to
+\infty} \ ||I(\cdot, t)||_{C([0, h(t)])}>0$, that is, spreading happens.
\end{thm}
\bpf We first consider the case that $R_0^F(0)(:=R_0^{DA}((-h_0, h_0),d_I, \alpha))>1$. In this case, we have  that the eigenvalue problem
\begin{eqnarray}
\left\{
\begin{array}{lll}
-d_I \psi_{xx}+\alpha \psi_x=\beta(x) \psi-\gamma(x)\psi+\lambda_0 \psi,\; &
x\in (-h_0,  h_0),  \\
\psi(x)=0, &x=\pm h_0
\end{array} \right.
\label{B2f}
\end{eqnarray}
admits a positive solution $\psi(r)$ with $||\psi||_{L^\infty}=1$, where $\lambda_0$ is the principle eigenvalue. It follows from Lemma 3.1
that $\lambda_0<0$.

We are going to construct a suitable lower solutions to
\eqref{a3} and we define
$$
\underline{I}(x,t)=
\delta \psi(x),\quad -h_0\leq x\leq h_0, \, t\geq 0,
$$
where $\delta$ is sufficiently small such that $0<\delta\leq \frac { N^*(-\lambda_0)}{\overline \beta}$
and $\delta \psi\leq I_0(x)$ in $[-h_0, h_0]$.

 Direct computations yield
\begin{eqnarray*}
& &\underline{I}_t-d_{I} \underline{I}_{xx}-(\beta(x)-\gamma(x))\underline{I}+\frac{\beta(x)}{N^*}\underline I^{2}\\
& &=\delta \psi(x)[\lambda_0+\frac{\beta(x)}{N^*}\delta \psi(x)]\\
& &\leq 0
\end{eqnarray*}
for all $t>0$ and $-h_0<x<h_0$.
Then we have
\begin{eqnarray*}
\left\{
\begin{array}{lll}
\underline{I}_t-d_{I}\underline{I}_{xx}+\alpha \underline{I}_x\leq (\beta(x)-\gamma(x))\underline{I}-\frac{\beta(x)}{N^*}\underline I^{2},\; &-h_0<x<h_0,\ t>0, \\
\underline{I}(\pm h_0,t)=0,\; &t>0, \\
0=h'_0\leq -\mu \underline{I}_{x}(h_0,t)=-\mu \delta \psi'(h_0),\  &t>0,\\
0=g'_0\geq -\mu \underline{I}_{x}(-h_0,t)=-\mu \delta \psi'(-h_0),\  &t>0,\\
\underline{I}(x,0)\leq I_{0}(x),\; &-h_0\leq x\leq h_0.
\end{array} \right.
\end{eqnarray*}
Hence we can apply Lemma 2.5 to conclude that  $I(x,t)\geq\underline I(x,t)$
in $[-h_0, h_0]\times [0,\infty)$. It follows that $\liminf_{t\to
+\infty} \ ||I(\cdot, t)||_{C([g(t), h(t)])}\geq \delta \psi(0)>0$ and therefore $h_\infty-g_\infty=+\infty$ by Lemma 4.1.

If $R_0^F(0)(:=R_0^{DA}((-h_0, h_0),d_I, \alpha))=1$. Then for any positive time $t_0$, we have $h(t_0)>h_0$, $g(t_0)<-h_0$ and $R_0^{DA}((g(t_0), h(t_0)), d_I, \alpha))>R_0^{DA}((-h_0, h_0), d_I, \alpha))=1$ by
the monotonicity in Theorem 3.4.
Replaced the initial time $0$ by the positive time $t_0$, we then have $h_\infty-g_\infty=+\infty$ as above.
 \epf

\begin{rmk} It follows from the above proof that spreading happens if there exists $t_0\geq 0$ such that $R_0^F(t_0)\geq 1$.
\end{rmk}

Next, we consider the long time behavior of the solution to problem \eqref{a3} when the spreading occurs.

\begin{thm} \label{van}If $h_\infty=-g_\infty=+\infty$, then the solution of
 free boundary problem \eqref{a3} satisfies $\lim_{t\to +\infty} \ I(x,t)=I^*(x)$
uniformly in any bounded subset of $(-\infty, \infty)$, where $I^*$ is the unique positive equilibrium of
the stationary problem:
\begin{eqnarray}
-d_I I^*_{xx}+\alpha I^*_x=(\beta(x)-\gamma(x))I^*-\frac{\beta(x)}{N^*} (I^*)^{2},\quad -\infty<x<\infty.
\label{odea3}
\end{eqnarray}
\end{thm}
\bpf We divide the proof in four parts.

(1) The existence and uniqueness of the stationary solution

It is easy to see that the comparison principle holds for the stationary problem with the solution in the sector $<0, N^*>:=\{I(x):\, 0\leq I(x)\leq  N^*,\,-\infty<x<\infty\}$.
Since that $h_\infty=-g_\infty=+\infty$, it follows from Remark 3.1 that there exists $t_0>0$ such that $R_0^F(t_0)=R_0^{DA}((g(t_0),h(t_0)),d_I,\alpha)>1$,
therefore, for any $l$ with $l\geq L_0:=\max\{-g(t_0), h(t_0)\}$, the problem
\begin{eqnarray}
-d_I I_{xx}+\alpha I_x=(\beta(x)-\gamma(x))I-\frac{\beta(x)}{N^*} (I)^{2},\ -l<x<l,\quad I(\pm l)=0
\label{odea31}
\end{eqnarray}
admits a unique positive solution $I_l(x)$. Using the comparison principle yields that as $l$ increases to infinity, $I_l$ increases to a positive solution of problem
(\ref{odea3}), which is referred as the minimal positive solution $\underline I^*$ of problem (\ref{odea3}).

On the other hand, any constant greater that $N^*$ is a upper solution of problem (\ref{odea3}), we then have the maximal positive solution $\overline I^*$ by the upper and lower solution method and the theory of monotone dynamical systems (  Corollary 3.6 in \cite{HS} or Theorem 5.1 in \cite{Pao}).

 The uniqueness ($\underline I^*=\overline I^*:=I^*$) of positive solution of problem
(\ref{odea3}) follows from the similar technique in \cite{DLiu} (Theorem 2.3) or \cite{Pao} (Theorem 5.3).

(2) The limit superior of the solution

It follows from the comparison principle that $I(x,t)\leq \overline I(x,t)$
for $(x,t)\in (-\infty, \infty)\times (0,\infty)$, where
$\overline I(x,t)$ solves
\begin{eqnarray}
\label{ode1}
\left\{
\begin{array}{lll}
&\overline I_t-d_I \overline I_{xx}+\alpha \overline I_x=(\beta(x)-\gamma(x))\overline I-\frac{\beta(x)}{N^*} (\overline I)^{2},&t>0\, -\infty<x<\infty, \\
&\overline I(x,0)=N^*,& -\infty<x<\infty.
\end{array} \right.
\end{eqnarray}
It is well known that $\overline I$ is monotone decreasing with respect to $t$ and $\lim_{t\to\infty}\overline
I(x,t)= I^*$ uniformly in any bounded subset of $(-\infty, \infty)$; therefore we deduce
\begin{equation}\label{123}
\limsup_{t\to +\infty} \ I(x,t)\leq I^*
\end{equation}
uniformly in any bounded subset of $(-\infty, \infty)$.

(3) The lower bound of the solution for a large time.

For part (1), we see that $I_{L_0}$ solves (\ref{odea31}) with $l$ replaced by $L_0$. Direct calculation shows that
we can choose $\delta$ sufficiently small such that $\delta I_{L_0}$ be a lower solution of the solution $I(x,t)$ in $[-L_0, L_0]\times [t_0, \infty)$.
We then have $I\geq \delta I_{L_0}$ in $[-L_0, L_0]\times [t_0, \infty)$,
which implies that the solution can not decay to zero, this result will be used in the next part.

(4) The limit inferior of the solution.

Since $h_\infty=-g_\infty=+\infty$,
 for any $L\geq L_0$, there exists $t_L>0$ such that $g(t)\leq -L$ and $h(t)\geq L$ for $t\geq t_L$.
We extend $I_{L_0}$ to $\phi_{L_0}(x)$ by defining $\phi_{L_0}(x):=I_{L_0}(x)$ for $-L_0\leq x\leq L_0$ and $\phi_{L_0}(x):=0$ for $x<-L_0$ or $x>L_0$.
Now for $L\geq L_0$, $I(x,t)$ satisfies
\begin{eqnarray}
\left\{
\begin{array}{lll}
I_t-d_I I_{xx}+\alpha I_x=(\beta(x)-\gamma(x)) I-\frac{\beta(x)}{N^*} I^{2},\; &  g(t)<x<h(t),\ t>t_L,  \\
I(x,t)=0, \quad & x=g(t)\, \textrm{or}\, x=h(t),\ t>t_L,\\
I(x,t_L)\geq \delta \phi_{L_0},
 & -L\leq x\leq L,
\end{array} \right.
\label{fs1}
\end{eqnarray}
therefore, we have $I(x,t)\geq w(x,t)$ in $[-L, L]\times [t_L, \infty)$,
where $w$ satisfies
\begin{eqnarray}
\left\{
\begin{array}{lll}
w_t-d_I w_{xx}+\alpha w_x=(\beta(x)-\gamma(x)) w-\frac{\beta(x)}{N^*} w^{2},\; &  -L<x<L,\ t>t_L,   \\
w(x,t)= 0, \quad & x=\pm L,\ t>t_L,\\
w(x,t_L)=\delta \phi_{L_0},  & -L\leq x\leq L.
\end{array} \right.
\label{fs11}
\end{eqnarray}
It follows from the upper and lower solution method
 and the theory of monotone dynamical systems ( \cite{HS} Corollary 3.6) that
$\lim_{t\to +\infty} \ w(x,t)\geq I_{L_0}(x)$ uniformly in
$[-L, L]$, where $I_{L_0}$ satisfies (\ref{odea31}) with $l$ replaced by $L_0$. Moreover,
By classical elliptic regularity theory and a diagonal procedure, it follows that  $I_{L}(x)$ converges
uniformly on any compact subset of $(-\infty, \infty)$ to $I^*(x)$.

Now for any given $[-M, M]$ with $M\geq L_0$, since that $I_{L}(x)\to I^*$ uniformly in $[-M, M]$, which is the compact subset of $(-\infty, \infty)$, as $L\to \infty$, we deduce that for any $\varepsilon >0$, there exists $L^*>L_0$ such that  $I_{L^*}(x)\geq  I^*-\varepsilon$ in $[-M, M]$.
As above, there is $t_{L^*}$ such that $[g(t), h(t)]\supseteq [-L^*, L^*]$ for $t\geq t_{L^*}$.
Therefore,
$$I(x,t)\geq w(x,t)\ \textrm{in}\ [-L^*, L^*]\times [t_{L^*}, \infty),$$ and
$$\lim_{t\to +\infty} \ w(x,t)\geq I_{L^*}(x) \textrm{in}\ [-L^*, L^*].$$
Using the fact that $I_{L^*}(x)(x))\geq I^*-\varepsilon$ in $[-M, M]$ gives
 $$\liminf_{t\to +\infty} \ I(x, t)\geq I^*(x)-\varepsilon \ \textrm{in}\ [-M, M].$$
 Since $\varepsilon>0$ is arbitrary, we then have $\liminf_{t\to +\infty} \ I(x, t)\geq I^*$  uniformly  in $[-M,M]$, which together with (\ref{123})
imply that $\lim_{t\to +\infty} \ I(x,t)=I^*$
uniformly in any bounded subset of $(-\infty, \infty)$.
\epf

\bigskip

Combing Lemma 4.1, Theorems 5.1 and 5.2, we immediately obtain the following
spreading-vanishing dichotomy:
\begin{thm}  Let $(I(x, t); g(t), h(t))$ be the solution of the  free boundary problem \eqref{a3}.
Then the following alternative holds:

Either
\begin{itemize}
\item[$(i)$] {\rm Spreading:} $h_\infty-g_\infty =+\infty$ and $\lim_{t\to
+\infty} \ I(x,t)=I^*$ uniformly for $x$ in any bounded set of $\R^1$, where $I^*$ is the unique positive solution of the stationary problem \eqref{odea3}; \end{itemize}

or
\begin{itemize}
\item[$(ii)$] {\rm Vanishing:} $h_\infty-g_\infty < \infty$, $R^{DA}_0((g_\infty, h_\infty), d_I, \alpha))\leq 1$ and $\lim_{t\to +\infty} \ ||I(
\cdot, t)||_{C([g(t),h(t)])}=0$.
\end{itemize}
\end{thm}
 \bpf
In fact, if $R^F_0(t_0)\geq 1$ for some $t_0\geq 0$, spreading happens by Theorem 5.1 or Remark 5.1.
Otherwise, $R^F_0(t)<1$ for any $t>0$, which means that $h_\infty-g_\infty<+\infty$,
vanishing happens.
 \epf

Theorem 4.2 shows if $R_0^F(0)<1$, vanishing happens for small expanding capability $\mu$ or small initial value of infected individuals, the next result
shows that spreading happens for large expanding capability and the proof will be omitted since it is an analogue of Lemma 3.7 in \cite{DL} or Lemma 2.8 in \cite{DG}.
 \begin{lem} Suppose that $R_0^F(0)<1$. Then $h_\infty-g_\infty=\infty$ if $\mu$ is sufficiently large.
\end{lem}

\begin{thm} (Sharp threshold)
Fixed $h_0$ and $I_0$. There exists $\mu^* \in [0, \infty)$
 such that spreading happens when $\mu> \mu^*$, and vanishing happens when $0<\mu\leq \mu^*$.
\end{thm}
\bpf
It follows from Theorem 5.1 that spreading always happens if $R_0^F(0)\geq 1$. Hence in this
case we have $\mu^*=0$.

For the remaining case $R_0^F(0)<1$. Define
$$\mu^*:=\sup \{\sigma_0: h_\infty (\mu)-g_\infty(\mu)<\infty \ \textrm{for}\ \mu\in (0,\mu_0]\}.$$
By Theorem 4.4, we
see that in this case vanishing happens for all small $\mu>0$, therefore, $\mu^*\in
(0, \infty]$. On the other hand, it follows from Lemma 5.4 that in
this case spreading happens for all big $\mu$. Therefore $\mu^*\in
(0, \infty)$, and spreading happens when $\mu> \mu^*$, vanishing happens when $0<\mu< \mu^*$ by Corollary 2.6.

We claim that vanishing happens when $\mu=\mu^*$. Otherwise $h_\infty-g_\infty=\infty$ for
$\mu=\mu^*$. Since $R_0^F(t)\to \sup_{r\in [0, \infty)}\frac {\beta (r)}{\gamma (r)}>1$ as $t\to \infty$,
 therefore there exists $T_0>0$ such that $R_0^F(T_0):=R_0^{DA}((g(T_0), h(T_0),d_I, \alpha)>1$. By
the continuous dependence of $(I,g, h)$ on its initial values, we can find
$\epsilon>0$ small so that the solution of (\ref{a3}) with $\mu=\mu^*-\epsilon$, denoted by
$(I_\epsilon,g_\epsilon, h_\epsilon)$ satisfies $R_0^{DA}((g_\epsilon(T_0), h_\epsilon(T_0)), d_I,\alpha)>1$. This implies that spreading
happens to $(I_\epsilon, g_\epsilon, h_\epsilon)$, which
contradicts the definition of $\mu^*$. The proof is complete.
 \epf

\section{Asymptotic spreading speeds}

To derive the asymptotic spreading speed, we first recall the known result for \eqref{a3} with $\alpha=0$,
see Corollary 3.7 in \cite{DG}.

\begin{thm} \label{spee}Let $(I; g, h)$ be the unique solution of \eqref{a3} with $\alpha=0$. If $h_\infty=-g_\infty=\infty$. Then
\[
\lim_{t\to\infty}\frac{-g(t)}{t}=\lim_{t\to\infty}\frac{h(t)}{t}=k_0,
\]
where $(k_0, q(x))$ is the unique positive solution of the problem
\begin{eqnarray}
\left\{
\begin{array}{lll}
-d_I q''+c_0 q'=q(a-bq), \quad &x>0,  \\
 q(0)=0,\ q(\infty)=a/b,\ q(x)>0, & x>0,\\
\;\;\; \mu q'(0)=k_0,&
\end{array} \right.
\label{m1}
\end{eqnarray}
\end{thm}
and $a=\beta_\infty-\gamma_\infty$, $b=\frac{\beta_\infty}{N^*}$.

\bigskip

Theorem \ref{spee} shows that if there is no advection, the asymptotic spreading speed of the left frontier and
that of the right frontier are the same when disease is spreading.

To address the change induced by an advection term, we first study the following problem:
\begin{equation}\label{prob-q-infty}
\left\{
\begin{array}{l}
d_I q'' -(k-\alpha)q' + q[a-b q]=0 \quad \mbox{ for } x\in (0,\infty),\\
q(0)=0, \ q(\infty)=\frac a{b}, \ q(z)>0\ \mbox{ for } x\in (0,\infty).
\end{array}
\right.
\end{equation}
 Usually, $q(z)$ is called a semi-wave with speed $k$. We will derive the rightward spreading speed by this semi-wave.
Consequently, for the leftward spreading speed, the corresponding semi-wave is governed by
\begin{equation}\label{prob-q2}
\left\{
\begin{array}{l}
d_I q'' -(k+\alpha)q' +q[a-b q]=0 \quad \mbox{ for } x\in (0,\infty),\\
q(0)=0, \ q(\infty)=\frac b{a}, \ q(z)>0\ \mbox{ for } x\in (0,\infty).
\end{array}
\right.
\end{equation}

We now present the properties of the semi-waves, see Propositions 2.2, 2.4 and 2.5 in \cite{GLL}.
\begin{prop}\label{prop-tw} The following conclusions hold.

\begin{itemize}
 \item[$(i)$] Problem \eqref{prob-q-infty} has exactly one solution
$(k,q) = (k^{*}_{r},q_{r}^{*})$ such that
\begin{equation}\label{c-q-r}
\mu (q_r^*)'(0)=k^*_r .
\end{equation}
Moreover, $k^{*}_{r}:=k^*_r(\alpha, d_I, a, b)\in (0,2\sqrt{ad_I}+\alpha)$;

 \item[$(ii)$] Problem \eqref{prob-q2} has exactly one solution
$(k,q) = (k^{*}_l,q_l^{*})$ such that
\begin{equation}\label{c-q-l}
\mu (q_l^*)'(0)=k^*_l .
\end{equation}
Moreover, $k^{*}_l:=k^*_l(\alpha, d_I, a,b)\in(0, 2\sqrt{ad_I}-\alpha)$;

\item[$(iii)$] $0 < k_l^* <k^* <k_r^* $, where $k^{*}$ is the speed in \eqref{prob-q-infty} (or \eqref{prob-q2} ) with $\alpha=0$;

\item[$(iv)$] $k_l^*$ and $k_r^*$
depend continuously on the parameter $a$, and are strictly increasing in $a$, that is,
for any $a>0$ and $a_1>a_2>0$, we have
$${k}_{r}^{*}(\alpha, d_I, a_1, b)>k_{r}^{*}(\alpha, d_I, a_2, b),\quad
\lim_{\varepsilon\to 0}{k}_{r}^{*}(\alpha, d_I, a+\varepsilon, b)=c_{r}^{*}(\alpha, d_I, a, b),$$
$${k}_{l}^{*}(\alpha, d_I, a_1, b)>k_{l}^{*}(\alpha, d_I, a_2, b),\quad
\lim_{\varepsilon\to 0}{k}_{l}^{*}(\alpha, d_I, a+\varepsilon, b)=c_{l}^{*}(\alpha, d_I, a, b);$$

\item[$(v)$] $k_l^*$ and $k_r^*$
depend continuously on the parameter $b$, and are strictly decreasing in $b$.

 \end{itemize}
\end{prop}

Next we give the spreading speeds when spreading happens.
\begin{thm}\label{spreadsp1}  If $h_\infty=-g_\infty=+\infty$, then
\begin{align*}
\lim_{t\to +\infty} \frac {h(t)}t=k^*_r(\alpha, d_I, a,b),\quad \lim_{t\to +\infty} \frac {-g(t)}t=k^*_l(\alpha, d_I, a,b),
\end{align*}
\label{m10}
where $a=\beta_\infty-\gamma_\infty$ and $b=\frac{\beta_\infty}{N^*}$.
\end{thm}

\bpf
 By assumption $(H)$, $\lim_{x\rightarrow\pm \infty}(\beta(x)-\gamma (x))=\beta_{\infty}-\gamma_{\infty}=a$,
$\lim_{x\rightarrow\pm \infty}\frac{\beta(x)}{N^*}=\frac{\beta_\infty}{N^*}=b$. Note that $\lim_{x\rightarrow\pm \infty}I^*(x)=b$.
Therefore, for any $\varepsilon>0$, there exists $L_1>0$ such that for $|x|\geq L_1$,
$$a-\varepsilon<\beta(x)-\gamma (x)<a+\varepsilon, \quad b-\varepsilon<\frac{\beta_\infty}{N^*}<b+\varepsilon.
$$
Owing to $\lim_{x\rightarrow\pm \infty} I^*(x)=\frac{(\beta_\infty-\gamma_\infty)N^*}{\beta_\infty}=\frac ab$, then for given $\varepsilon$,
there exists $L_2>L_1$ such that for $|x|\geq L_2$,
$$\frac a b-\varepsilon <I^*(x)<\frac ab+\varepsilon.$$

Using the comparison principle and following the proof of Theorem 3.6 in \cite{DG}, we can get

\begin{equation}\label{est-lower}
\liminf\limits_{t\to \infty}   \frac{h(t)}{t} \geq {k}^*_r(\alpha, d_I, a-\varepsilon, b+\varepsilon) .
\end{equation}

\begin{equation}\label{est-upper}
\limsup\limits_{t\to \infty}   \frac{h(t)}{t} \leq k_r^*(\alpha, d_I, a+\varepsilon, b-\varepsilon).
\end{equation}
Letting $\varepsilon \to 0$ give that
\begin{equation}\label{est-upper1}
\limsup\limits_{t\to \infty}   \frac{h(t)}{t} \leq k_r^*(\alpha, d_I, a, b)
\end{equation}
by Proposition \ref{prop-tw} $(iv)$ and $(v)$.

The leftward spreading speed can be discussed similarly.

\epf

\section{\bf Numerical illustration and discussion}

In this section, we first carry out numerical simulations to illustrate the theoretical results given above.
Because the boundary is unknown, it is a little difficult to present the numerical solution compared to the problem in fixed boundary.
We use an implicit scheme as in \cite{Ra2009} and then obtain a nonlinear
system of algebraic equations, which was solved with Newton-Raphson method.

Let us fix some coefficients and functions. Assume that
$$N^*=2,\ d_I=4,\ h_0=1,\ I_0(x)=\cos (\frac \pi 2 x),$$
$$\beta(x)=4+\frac 2{1+x^2}\sin x,\ \gamma(x)=1+\frac 1{1+x^2}\cos x,$$
then the asymptotic behaviors of the solution to problem $(\ref{a3})$
are shown by choosing different advection rate $\alpha $ and expanding capability $\mu$.

 \begin{exm}
 Fix big expanding capability $\mu=6$, and choose $\alpha =1.5$ and $\alpha=-1.5$, it is easy to see from Figure \ref{t1} that the free boundaries $x=g(t)$ and $x=h(t)$ increase fast, and the solution $I$ stabilizes to a positive equilibrium. Moreover, owing to the advection, the right boundary goes faster that the left one in the left graph for $\alpha =1.5$. Contrarily, in the right graph, $\alpha=-1.5$ and the left boundary goes faster.
\end{exm}

\begin{figure}
\centerline{ \includegraphics[width=2.8in]{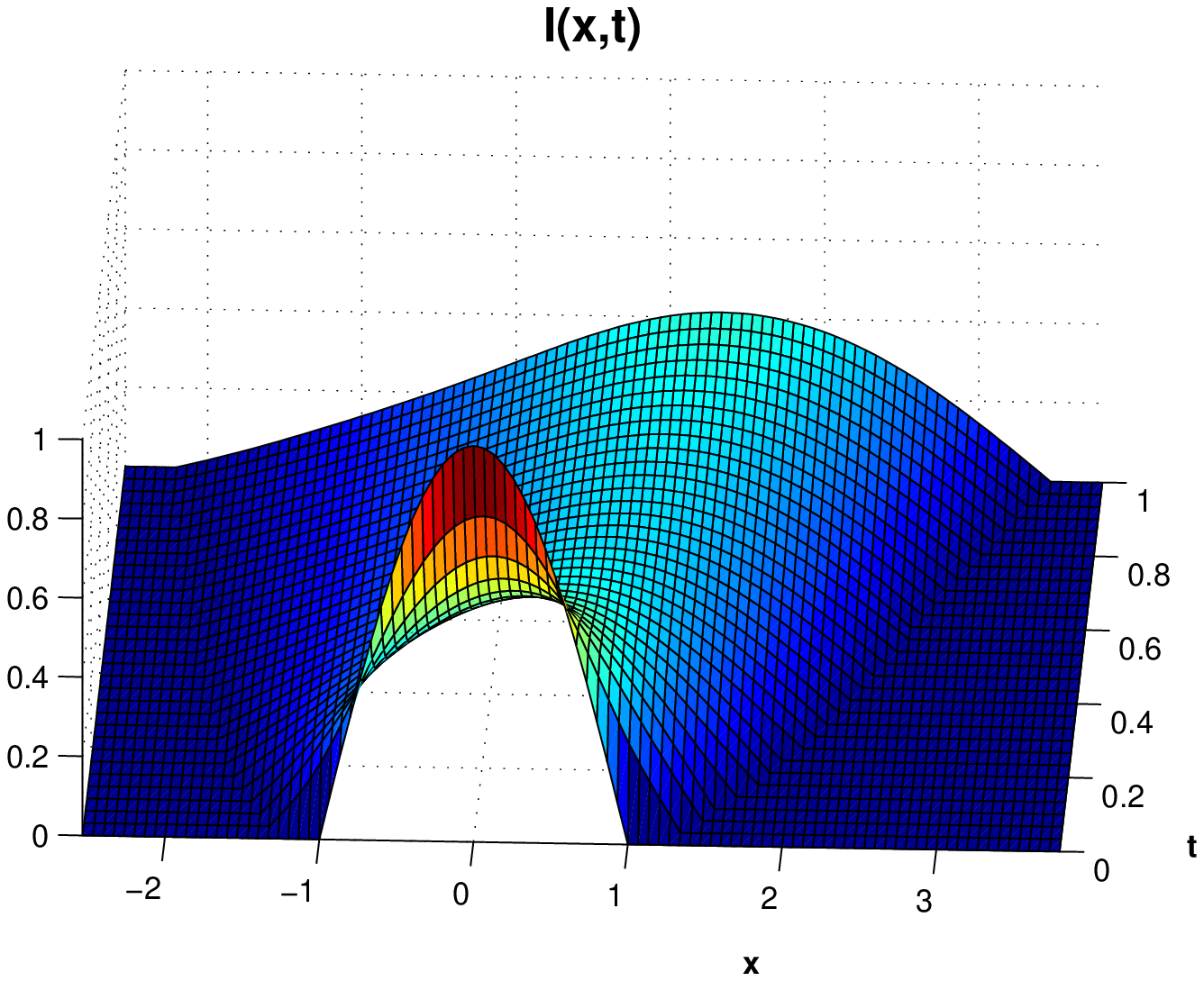}
\includegraphics[width=2.8in]{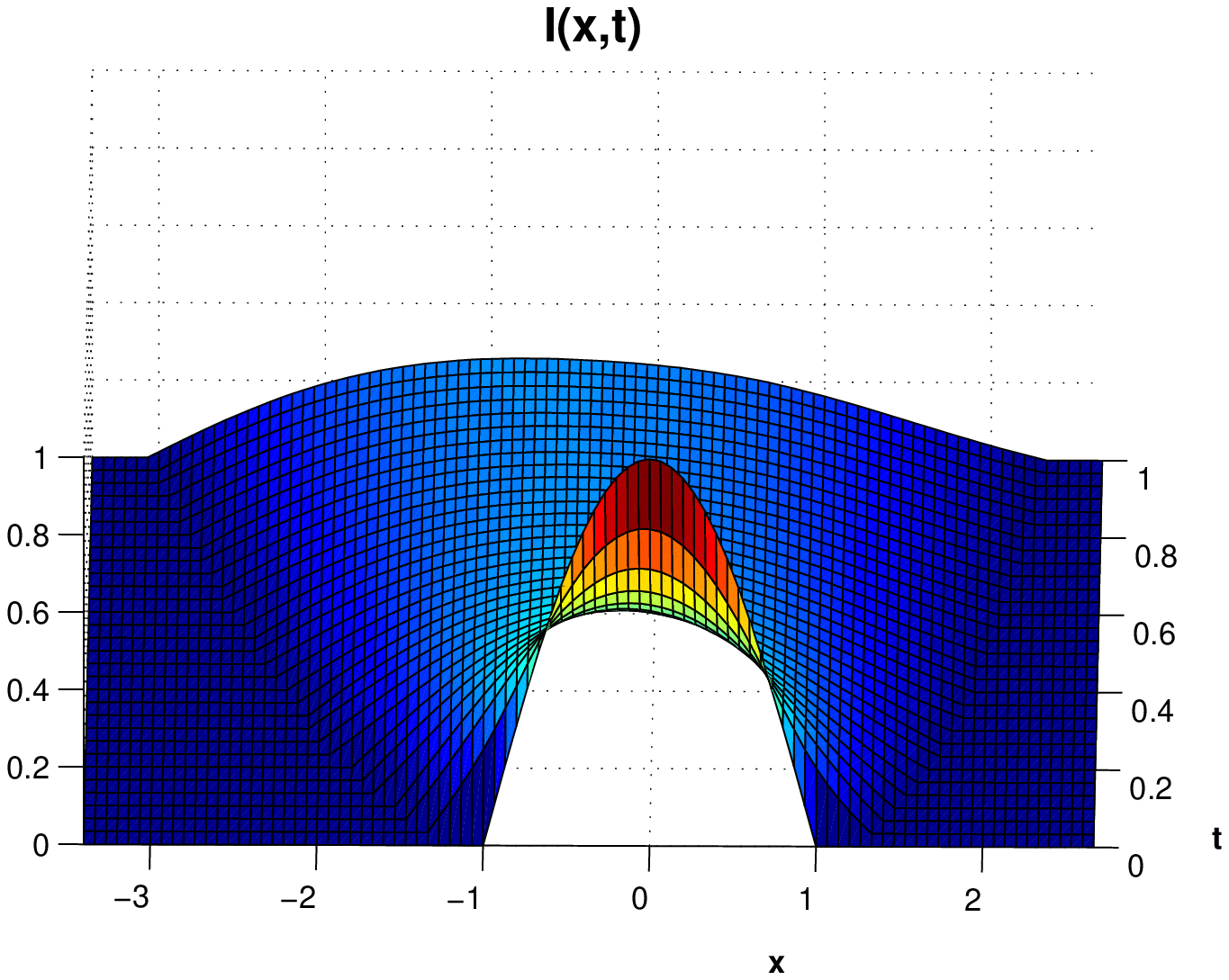}}
  \caption{ $\mu=6$, $\alpha =1.5$ for the left and $\alpha=-1.5$ for the right. The solution $I$ in the left graph turns right and stabilizes to a positive equilibrium, while in the right graph, $I$ turns left and stabilizes to a positive equilibrium.}\label{t1}
\end{figure}

 \begin{exm}
Fix small $\mu=1$, and choose $\alpha =1.5$ and $\alpha=-1.5$, compared the free boundary in Figure \ref{t2} with that in Figure \ref{t1}, the free boundaries $x=h(t)$
and $g(t)$  in Figure \ref{t2} increase slower than that in Figure \ref{t1}. Moreover, the solution $I$ decays to zero quickly.
\end{exm}

\begin{figure}
\centerline{ \includegraphics[width=2.8in]{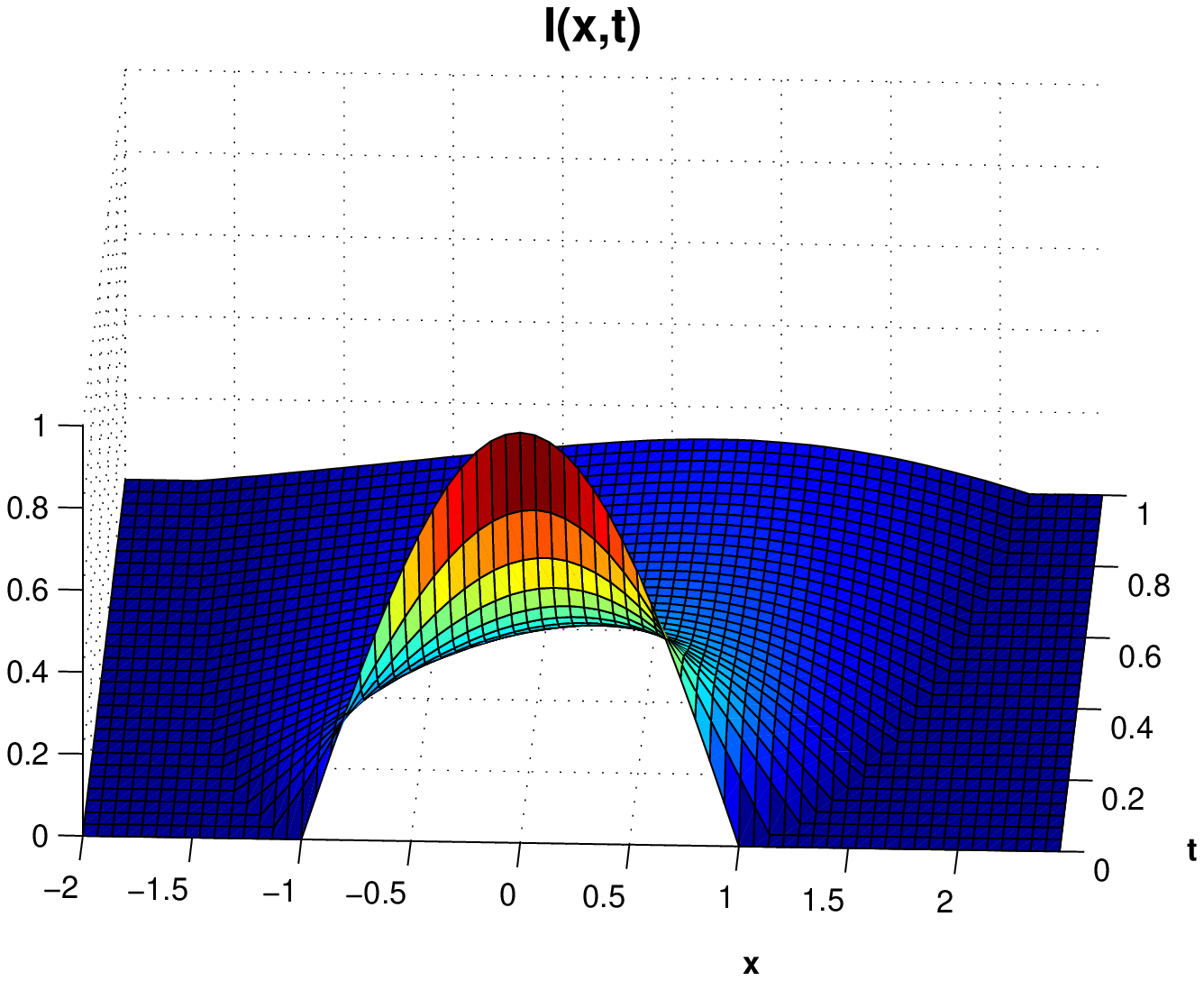}
\includegraphics[width=2.8in]{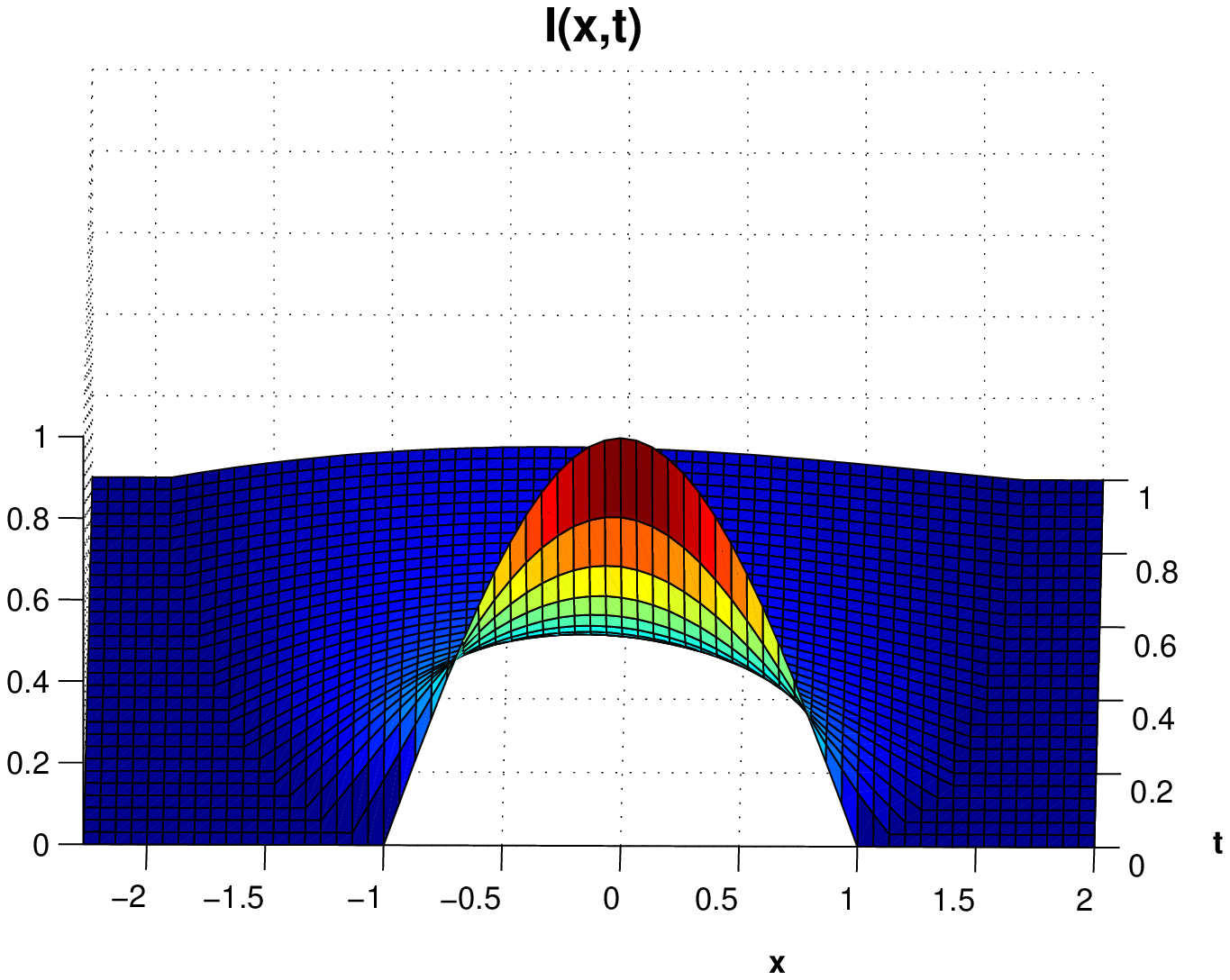}  }
  \caption{$\mu=1$, $\alpha =1.5$ for the left and $\alpha=-1.5$ for the right. The solution $I$ decays to zero quickly and the free boundaries increase slowly.}\label{t2}
\end{figure}

In this paper, we have considered a simplified spatial SIS epidemic model
describing the spatial transmission of diseases and examined the dynamical behavior of the population $I$
with spreading fronts $x=h(t)$ and $x=g(t)$ defined by \eqref{a3}.  We have
obtained some analytic results about the asymptotic properties of the spatial spread of infectious diseases.

The basic reproduction numbers $R_0^{DA}$ and $R_0^F(t)$ are introduced for the diffusion-reaction-advection system with Dirichlet boundary condition
and the system with the free boundary, respectively.  It is proved that if $R_0^F(t_0)\geq 1$ for some $t_0\geq 0$, spreading always
happens or the disease will become endemic (Theorem 5.1 and Remark 5.1).
If $R_0^F(0)<1$, vanishing of the spreading of the disease happens provided that the initial value of the infected individuals $I_0$ is sufficiently
small (Theorem 4.3) or the expanding capability is small (Theorem 4.4), while spreading happens provided that the expanding capability is large (Lemma 5.4).

One of the main contributions of this work is the development and analysis of the basic reproduction
numbers. We now have four basic reproduction
numbers: $R_0$ used for the ODE system, $R_0^N$ (see (\ref{AAq})) defined for the diffusive system with
Neumann boundary condition,  $R_0^{DA}$ and $R_0^F(t)$ defined in this paper.
They are all closely related, $R_0$ is actually equals to $R_0^N$ with $\Omega$ replaced by the whole space ${\mathbb R}^n$ if all coefficients are constant, $R_0^F(t)$ is $R_0^{DA}$ with
$\Omega$ replaced by the changing interval $(g(t),h(t))$. However, they are different,  $R_0$, $R_0^N(t)$ and $R_0^{DA}$ are all constants, while
$R_0^F(t)$ depends on time $t$, the temporal dependence of the basic reproduction number is a intrinsic characteristic of the spreading over a changing domain.
It follows from the definition of $R_0^{DA}$ (or $R_0^F(0)$)  that fast diffusion
 and small initial infected size are in favor of the disease to vanish, or prevention and control, the latter implies that early control is better to prevent the outbreak of the disease to spread over larger area.

 Another consideration of this work is the impact of advection on the left and right free boundaries. Because of wind direction, human activities and the migration of birds, etc., disease prefer to move towards one direction. Introduction of the small advection in this paper reveals the different asymptotic spreading speeds,
 big advection, we believe, will causes more complex dynamical behaviors. We keep it as a future work when use West Nile virus as a concrete example.

\end{document}